\documentclass[a4paper]{article}

\usepackage[T1]{fontenc}
\usepackage[utf8]{inputenc}
\usepackage[ngerman]{babel}
\usepackage{latexsym,amssymb,amsfonts,amsmath,mathabx}
\usepackage{csquotes}
\usepackage{chngcntr}
\usepackage{lmodern}
\usepackage[absolute]{textpos}
\usepackage{float}
\usepackage{graphicx}
\usepackage{hyperref}
\usepackage{orcidlink}
\usepackage{authblk}

\usepackage[style=numeric, backend=biber]{biblatex}
\begin{document}

\title{Über Nichts}

\author[András Bátkai]{András Bátkai \orcidlink{0000-0002-9209-2779}}

\affil{Institut für Sekundarbildung und Fachdidaktik, Pädagogische Hochschule Vorarlberg, Liechtensteiner Str. 33-37, A-6800, Feldkirch, Österreich, \href{mailto:andras.batkai@ph-vorarlberg.ac.at}{\texttt{andras.batkai@ph-vorarlberg.ac.at}}}

\maketitle

\begin{abstract}
In diesem Beitrag wird das Konzept des \glqq Nichts\grqq{} an der Schnittstelle von Philosophie und Mathematik untersucht. Ausgehend von philosophischen Deutungen – insbesondere Platons relationalem Nichts und Kants \textit{ens imaginarium} – wird gezeigt, wie die Mathematik formale Absenz in existierende, operative Strukturen übersetzt. Im Zentrum stehen dabei die Zahl Null (0) und die leere Menge ($\emptyset$). Die historische und kognitive Entwicklung der Null verdeutlicht ihren Wandel vom blo\ss{}en Platzhalter zur ontologischen Symmetrieachse und zum produktiven Werkzeug der Algebra. Parallel dazu erweist sich die leere Menge nicht als absoluter Mangel, sondern als axiomatisch zwingendes Fundament und kategorientheoretisches Initialobjekt. Darüber hinaus wird die erkenntnistheoretische Dimension des Nichts beleuchtet: Vom sokratischen Nichtwissen als schöpferischem Vakuum spannt sich der Bogen über Nikolaus von Kues bis hin zur reinen Geometrie Bólyais und der Mengenlehre Cantors. Das Nichts offenbart sich letztlich nicht als Abgrund der Bedeutungslosigkeit, sondern als jener unentbehrliche Raum, der die Konstruktion der Mathematik als freie Schöpfung des menschlichen Geistes und als \glqq Wissenschaft von Mustern\grqq{} überhaupt erst ermöglicht.
\end{abstract}


\setcounter{footnote}{0}
\setcounter{section}{0}
\setcounter{figure}{0}
\setcounter{equation}{0}

\section{Einleitung: Nichts und Mathematik}

Das \glqq Nichts\grqq{} wurde in vielerlei Hinsicht sowohl mathematisch als auch philosophisch in den letzten Jahrtausenden intensiv untersucht. Fragen, wie man Leere, Nichts, Absenz formalisiert, damit umgeht, einordnet, haben Menschen lange beschäftigt und zu vielen Kontroversen geführt. Der Begriff des Nichts steht im Zentrum der metaphysischen Reflexion und fordert die Mathematik auf ihren fundamentalsten Ebenen heraus. Bevor wir unsere Reise um das mathematische Nichts herum beginnen, sollten wir kurz klären, worum es hier geht. In Metzlers Lexikon der Philosophie \cite{metzler} lesen wir unter dem Schlagwort \glqq Nichts\grqq{} folgende Verwendungen:
\begin{enumerate}
\item Nichts als absolutes Nichts (Nihil).
\item Nichts als Privation, Änderung, Mangel.
\end{enumerate}

 Diese noch rein binäre Unterscheidung lässt sich jedoch durch Immanuel Kants \glqq Tafel vom Begriffe des Nichts\grqq{} in der \textit{Kritik der reinen Vernunft} (A 290/B 347) wesentlich präziser fassen und systematisch ordnen \cite{kant_krv}. Kant unterscheidet dort vier Arten des Nichts: 
\begin{enumerate}
    \item das \textit{ens rationis} (leerer Begriff ohne Gegenstand), 
    \item das \textit{nihil privativum} (leerer Gegenstand eines Begriffs, d. h. die Negation als Mangel einer Realität wie Schatten oder Kälte), 
    \item das \textit{ens imaginarium} (leere Anschauung ohne Gegenstand, wie der reine Raum oder die reine Zeit) und 
    \item das \textit{nihil negativum} (der Gegenstand eines sich selbst widersprechenden Begriffs, das logisch Unmögliche).
\end{enumerate}
Während das \textit{nihil negativum} die Grenze der logischen Denkbarkeit markiert, bieten insbesondere das \textit{nihil privativum} und das \textit{ens imaginarium} die philosophische Basis, um die mathematische Null und die leere Menge nicht als absolutes Nichts, sondern als formale Strukturen und Bedingungen der Erkenntnis zu begreifen.

Die erste Verwendung, das philosophische Nichts (Nihil), wird klassisch als radikale Negation des Seins verstanden, wie es schon die antiken Denker Parmenides und Demokrit behandelten. Diese ersten Diskussionen warfen ein Licht auf das Paradoxon zwischen dem Nichts und der Existenz des Nichts. Die zweite Verwendung des Nichts als Negation und Änderung ist untrennbar mit Platons Namen verbunden, insbesondere mit seinem Spätwerk Sophistes \cite[257B--258E]{2007platon}. In diesem Dialog untersuchen der Fremde aus Elea und der junge Mathematiker Theaitetos das Wesen des Sophisten und sto\ss{}en dabei auf das Problem der Definition des Nichtseienden. Platon positioniert das Nichts hier nicht als absolute Abwesenheit, sondern lässt den Fremden argumentieren, dass wir mit dem Begriff des \glqq Nichtseienden\grqq{} nicht das blo\ss{}e Gegenteil des Seienden meinen. Der Gast aus Elea stellt fest:

\begin{quote}
Wenn wir \glq Nichtseiendes\grq{} sagen, so meinen wir nicht, wie es scheint, ein entgegengesetztes des Seienden, sondern nur ein verschiedenes. (257B)
\end{quote}

Diese Neudeutung ist fundamental: Das Nichts wird zu einem aktiven, konstituierenden Prinzip der Andersheit. Im Verlauf des Gesprächs wird diese relative Existenz des Nichtseienden weiter gefestigt. Der Fremde führt Theaitetos zu der Einsicht, dass das Nichtseiende unbestritten sein eigenes Wesen besitzt – analog dazu, wie das \glqq Nicht-Gro\ss{}e\grqq{} oder \glqq Nicht-Schöne\grqq{} genauso real ist wie das Gro\ss{}e oder Schöne selbst. Er fragt Theaitetos:

\begin{quote}
Und darf man schon sagen, dass das Nichtseiende unbestritten sein eigenes Wesen hat, und so wie das Gro\ss{}e gro\ss{} und das Schöne schön war, und das Nichtgro\ss{}e und Nichtschöne nichtgro\ss{} und nichtschön, ebenso auch das Nichtseiende war und ist nichtseiend und mit einzuordnen als ein Begriff unter das Viele das ist? (258B-C)
\end{quote}

Durch diese unterscheidende Funktion des Nichts werden die Vielfalt der Existenz und die Verflechtung von Ideen überhaupt erst erklärbar. Am Ende ihrer Untersuchung können die Dialogpartner somit folgern, dass das Nichtseiende nicht nur existiert, sondern als die Natur der Verschiedenheit zu begreifen ist. Der Gast fasst dies abschlie\ss{}end zusammen:

\begin{quote}
Wir aber haben nicht nur gezeigt, dass das Nichtseiende ist, sondern haben auch den Begriff des Nichtseienden aufgezeigt. Denn nachdem wir gezeigt, dass die Verschiedenheit ist, und dass sie verteilt ist unter alles Seiende gegeneinander, so haben wir von jedem Seienden entgegengesetzte Begriffe derselben zu sagen gewagt, dass eben dies in Wahrheit das Nichtseiende sei. (258D-E)
\end{quote}

Die Mathematik hat naturgemä\ss{} mehr Kontaktpunkte mit dieser zweiten, relationalen Verwendung des Nichts. Sie formalisiert die Absenz und weist ihr eine aktive, systematisch notwendige Existenz in Form der Zahl Null (0) und der leeren Menge ($\emptyset$) zu.  Im Folgenden untersuchen wir zunächst diese beiden fundamentalen mathematischen Modellierungen und zeigen auf, wie das Nichts in der mathematischen Praxis nicht als blo\ss{}er Mangel, sondern als nützliches, produktives und operatives Werkzeug fungiert. Abschlie\ss{}end beleuchten wir die erkenntnistheoretische Dimension des Nichts: Ausgehend vom sokratischen Nichtwissen als schöpferischem Nullpunkt präsentieren wir anhand der Philosophie von Nikolaus von Kues sowie der späteren Arbeiten von János Bólyai, Georg Cantor, George Spencer Brown und Imre Lakatos den Gedankengang, wie sich die Mathematik als freie Schöpfung des menschlichen Geistes ex nihilo entfaltet und letztlich als reine \glqq Wissenschaft von Mustern\grqq{} begreifen lässt.

Der Autor bedankt sich beim anonymen Gutachter für das aufmerksame Lektorat und die scharfsinnigen Verbesserungsvorschläge, die ma\ss{}geblich zur Qualität dieser Arbeit beigetragen haben.

\section{Nichts und die Null}

Die Geschichte der Zahl Null ist die Geschichte ihrer fortschreitenden Konzeptualisierung vom blo\ss{}en Fehlen zum aktiven Fundament. Die Mathematikdidaktik formuliert folgende Aspekte und Deutungen der Null \cite{fahse2014}:

\begin{itemize}
\item Kardinalzahlaspekt (Mengenkonzept): Null steht für \glqq nichts da\grqq{} oder \glqq gar keine\grqq. Zum Beispiel: \glqq Wie viele Äpfel hast du?\grqq{} – \glqq Null\grqq.
\item Operatoraspekt (\glqq Wie oft?\grqq): Null als \glqq gar keine\grqq{} Wiederholung, z.B. \glqq 0 mal 7 ist Null\grqq{} bedeutet, dass die Aktion \glqq Lege eine Siebenerreihe hin\grqq{} nicht ausgeführt wurde.
\item Kodierungsaspekt (Stellenwertsystem): Null als Platzhalter, der das Nicht-Vorhandensein einer Stelle in einer Zahl anzeigt, z.B. in der Zahl 102.
\item Ma\ss{}zahlaspekt: Null als Ma\ss{} für eine Grö\ss{}e, wie $0$ Euro auf ein Konto, $0$ Meter Entfernung.
\item Referenzaspekt: Ein willkürlich gewählter Punkt zur Orientierung, wie z.B. eine Temperatur von  $0,00^{\circ}\text{C}$, zero milestone in Washington D.C., oder der Nullmeridian durch Greenwich, England.
\item Rechenaspekt: Null als Zahl, mit der gerechnet wird (Addition, Subtraktion, Multiplikation), die aber auch Stolpersteine birgt.
\end{itemize}
Manchmal kommt der Ordinalzahlaspekt dazu, wie nullte Stunde in der Schule, oder in Gebäuden Erdgeschoss als nulltes Geschoss. Eine besondere Schwierigkeit sehen wir in der Volksschule, wo der Kodierungsaspekt viel früher auftritt als der Kardinalzahlaspekt.

Jenseits der funktionalen Aspekte stellt sich die Frage nach der neurobiologischen Basis der Null. Die kognitionswissenschaftliche Forschung, insbesondere die Arbeiten von Nieder \cite{nieder2016}, zeigt, dass das menschliche Gehirn eine Brücke von der Wahrnehmung der blo\ss{}en Absenz hin zur abstrakten Zahl schlagen muss. Während die \textit{Absence Perception} (Wahrnehmung von Abwesenheit) ein biologisches Grundvermögen ist, erfordert die Zahl Null eine neuronale Kategorisierung des Nichts als \glqq Etwas\grqq{}. Jüngste Studien \cite{kutter2024} belegen, dass spezifische Neuronen im medialen Temporallappen die Null als kleinsten Wert auf einer mentalen Zahlenlinie kodieren. Die Null ist somit im Gehirn nicht durch ein Schweigen der Neuronen repräsentiert, sondern durch ein aktives Feuern, das eine quantitative Kategorie für die Leere schafft.

Diese kognitive Fähigkeit zur Abstraktion spiegelt sich historisch in der Entwicklung der ersten Schriftsysteme wider: Als die erste Erscheinung von Null als Ziffer in dem frühesten Stellenwertsystem der Altbabylonier könnte ein Zeichen (der \glqq Doppelnagel\grqq) oder eventuell ein leerer Raum \cite[Seite 20]{neugebauer1969} zwischen den Ziffern lediglich als Platzhalter innerhalb einer Zahl, um die relative Position der Ziffern zueinander anzuzeigen, gewertet werden. Die Verwendung dieser Platzhalterbezeichnungen ist eher die Ausnahme, oft haben die Schreiber keinen Bedarf gesehen. Als Ziffer im Stellenwertsystem wurde die Null (der \glqq Doppelnagel\grqq) als Zeichen ca. 400 v. Chr. durchgehend in astronomischen Texten \cite{neugebauer1941} verwendet. Hier erscheint die Null das erste Mal als End- (oder Anfangs-) Null, um die absolute Positionierung der Ziffern und somit die Grö\ss{}enordnung der Zahlen anzuzeigen. Dieses Symbol wurde  früher auch für die Satztrennung benutzt.

Im altbabylonischen Zahlensystem gab es noch kein Symbol für die Null als Anzahl. Das Ergebnis von Subtraktionen wie 20-20 war einfach nicht geschrieben \cite[Seite 254]{friberg2016} oder als \glqq fehlend\grqq{} bezeichnet \cite[Seite 293]{hoyrup2013} worden. Auch in späteren mathematischen Texten aus Susa wurde das umschrieben (\glqq du wei\ss{}t ja\grqq) oder bei Getreideverteilungen als \glqq Das Korn ist ausgegangen\grqq{} vermerkt – eine deutliche Indikation, dass der Begriff des Nichts noch nicht in Form einer Zahl gefasst war.

Die entscheidende Wende vollzogen die Inder \cite[Seite 438]{ifrah2000} vor etwa 1500 Jahren. Sie bezeichneten mit dem Sanskritwort śūnya explizit die Zahl \glqq Null\grqq{} und zugleich \glqq die Leere\grqq. Ihnen verdanken wir die dreifache Erfindung der Null, wie wir sie heute kennen: Als Ziffer in Stellenwertsystemen (Platzhalter und Andeutung der Zehnerpotenzen-Bündelung), als Bezeichnung einer Menge und als Zahl selbst (als Ergebnis der Subtraktion).

Irgendwann zwischen 3. Jahrhundert v. Chr. und 7. Jahrhundert n. Chr. wurden in Indien negative Zahlen verwendet, um Schulden darzustellen, und positive Zahlen, um Vermögenswerte zu kennzeichnen. Auf diese Weise konnte ein Mangel oder eine Schuld als etwas Reales konzeptualisiert werden, ebenso wie der Besitz von etwas, das – ob materiell oder nicht – real und zählbar ist. Eine wirtschaftlich und fiskalisch fortgeschrittene Kultur war fast sicher notwendig für die Entwicklung negativer Zahlen als Zahlen, die den natürlichen (oder positiven) Zahlen gleichgestellt sind. Dies erforderte jedoch einen bedeutenden konzeptuellen Schritt über die Vorstellung hinaus, dass Zahlen etwas sinnlich materiell Gegenwärtiges oder potenziell Gegenwärtiges repräsentieren. Denn quantifizierter Besitz kann individuell durch die Sinne überprüft werden, wohingegen finanzielle Vermögenswerte oder Schulden nur durch Bezugnahme auf soziale Vereinbarungen und Dokumente, die sie aufzeichnen, bestimmt werden können. Daher stellt die Verwendung von Schulden als negativen Zahlen einen signifikanten konzeptuellen Fortschritt in der Entwicklung von Zahlensystemen dar.

Wir kennen heute diese Errungenschaften aus den Werken von Brahmagupta (598–665 n. Chr.), der analog zu Euklid \cite{Chapter17NoughtMatterstheHistoryandPhilosophyofZero} eine systematisierende Rolle einnimmt. Er war bereits mit negativen Zahlen vertraut, als er die Eigenschaften der Null als eigenständige Zahl formulierte. Er definierte sie als die Zahl, die man erhält, wenn man eine Zahl von sich selbst subtrahiert. Damit vollzog (oder nutzte) er einen weiteren Abstraktionsschritt von Zahlen, der das Ergebnis des Zählens materiell vorhandener Objekte ist. Über die Konzepte von Schuld und Leere hinauszugehen und die Null als Zahl auf einer Stufe mit anderen Zahlen zu konzeptualisieren, ist immer noch ein riesiger Gedankenschritt.

Brahmagupta akzeptierte, dass sowohl $2-3$ als auch $2-2$ legitime Operationen \cite[Seite 439]{ifrah2000} sind und dass jede eine erkennbare Zahl als Lösung liefert. Er schuf eine ausgefeilte Theorie, die all diese Komponenten umfasste und damit mehr oder weniger den modernen Bereich der ganzen Zahlen ($\mathbb{Z}$) in seiner heutigen Form etablierte. Brahmagupta arbeitete auch in der Algebra und betrachtete Null- und Negativlösungen als akzeptabel und legitim. Indem er das Wissen über positive und negative Zahlen sowie die Null und all ihre Beziehungen zusammenführte und erweiterte, vollbrachte Brahmagupta eine gro\ss{}e Synthese und einen riesigen Sprung \cite{Chapter17NoughtMatterstheHistoryandPhilosophyofZero} nach vorn.

Diese mathematische Konzeptualisierung negativer Werte findet ihre philosophische Entsprechung in Kants nihil privativum. Wie Kant in seiner Schrift \textit{Versuch den Begriff der negativen Grö\ss{}en in die Weltweisheit einzuführen} (1763) \cite{kant2013} detailliert ausführt, ist das Nichts in diesem Kontext nicht als absolute Abwesenheit zu verstehen, sondern als reale Entgegensetzung (Realrepugnanz). Ein prägnantes physikalisches Beispiel für eine solche reale Entgegensetzung – das über das von sozialen Konventionen abhängige Konzept von Besitz und Schuld hinausgeht – ist das Verhältnis von positiver und negativer elektrischer Ladung. Die Null fungiert in diesem System nicht als blo\ss{}es Nichts, sondern als der kritische Grenzpunkt und der Zustand des Gleichgewichts, an dem sich zwei entgegengesetzte Realitäten exakt aufheben.

Mit der Etablierung des Bereichs der ganzen Zahlen ($\mathbb{Z}$) wandelt sich die Rolle der Null von einer blo\ss{}en Leerstelle hin zu einer ontologischen Symmetrieachse. Sie fungiert als der \glqq Umsprungpunkt\grqq{} zwischen Sein und einem spiegelbildlichen Gegen-Sein. In einer Philosophie des Nichts lässt sich die Null hier als jener Punkt definieren, an dem sich entgegengesetzte Qualitäten gegenseitig aufheben — ein Konzept, das in der modernen Physik seine Entsprechung in der Annihilation von Materie und Antimaterie findet. Diese Symmetrie-Konzeption \cite{kaplan2000} begreift die Null nicht als Vakuum, sondern als einen Zustand perfekten Gleichgewichts. Auf der Zahlengeraden ist sie die neutrale Mitte, die positive Werte und negative Schulden voneinander scheidet und gleichzeitig verbindet. Dieser Aspekt wird uns in Abschnitt \ref{sec:nutz} erneut beschäftigen.

Das Wort für Null behält bis heute in vielen Sprachen (zero) deutliche Spuren seiner hinduistischen und arabischen Wurzeln. Als die Araber die hindu-arabischen Ziffern übernahmen, übernahmen sie auch die Null, deren Begriff \glqq sunya\grqq{} sie in sifr umwandelten. Einige westliche Gelehrte wandelten sifr in ein lateinisch klingendes Wort um, nämlich zephirus, welches die Wurzel des Wortes zero ist. Andere westliche Mathematiker nannten die Null cifra, woraus das Wort cipher (Ziffer/Chiffre) wurde. Aufgrund der fundamentalen Bedeutung der Null für das neu etablierte Stellenwertsystem verallgemeinerte sich dieser Begriff im Laufe der Zeit, sodass das ursprünglich nur für die Null stehende Wort (cifra / Ziffer) schlie\ss{}lich als Bezeichnung für alle Zahlzeichen übernommen wurde\footnote{Über die historische, philosophische und kulturelle Bedeutung von Null verweisen wir auf den Sammelband \cite{TheOriginandSignificanceofZero} und auf das Buch \cite{rotman1987}.}.

Ein bemerkenswerter, von der eurasischen Entwicklung isolierter Entwurf der Null findet sich in der Mathematik der Maya (ca. 300–900 n. Chr.). Während die indische sunya primär aus dem Konzept der Leere hervorging, ist die Maya-Null (nik) eng mit dem Konzept der Vervollständigung und zyklischen Zeit verknüpft \cite[Chapter 30]{ifrah2000}. In ihrem vigesimalen Stellenwertsystem wurde die Null oft durch ein muschelartiges Symbol oder ein stilisiertes Schneckenhaus dargestellt. Obwohl in der Praxis als Platzhalter benutzt, philosophisch markiert die Maya-Null keinen blo\ss{}en Mangel, sondern den Abschluss einer Einheit, bevor ein neuer Zyklus beginnt \cite[Chapter 3]{sharer2006}. In den Inschriften der \glq Langen Zählung\grq{} wird deutlich, dass die Null notwendig war, um die Kontinuität der Zeitstrukturen zu wahren – ein Aspekt, den Chrisomalis \cite{chrisomalis2003} als strukturelle Notwendigkeit der Abwesenheit beschreibt. Damit bietet die Maya-Null eine komplementäre Sichtweise zum westlichen Verständnis: Das \glq Nichts\grq{} ist hier nicht der Abgrund, sondern das Fundament der zyklischen Ordnung.

Die Null ist in der modernen Arithmetik keineswegs nichtig, sondern trägt eine besondere Rolle:
\begin{itemize}
\item Neutrales Element: Bei der Addition verkörpert sie das Statische, da sie ein Element unverändert lässt ($a + 0 = a$).
\item Absorbierendes Element: Bei der Multiplikation wirkt sie absorbierend ($a \cdot 0 = 0$).
\end{itemize}
Darüber hinaus zeichnet sich die Null durch ihre singulären strukturellen Eigenschaften aus. 

Ein kritischer Grenzpunkt der mathematischen Logik manifestiert sich in der Division durch Null. Während die Addition des Nichts das System stabil hält, führt die Division durch die Null zu einem strukturellen Kollaps der Arithmetik. Philosophisch betrachtet markiert diese Unmöglichkeit eine logische Singularität: Da die Division als wiederholte Subtraktion verstanden werden kann, würde der Versuch, eine Menge durch das Nichts zu teilen, einen unendlichen Prozess auslösen. 

In der Terminologie Kants lie\ss{}e sich dieser Zusammenbruch der Arithmetik als Begegnung mit dem nihil negativum deuten. Während die Zahl Null für sich genommen ein widerspruchsfreier Begriff ist, erweist sich die Operation der Division durch sie als ein logisches \glqq Unding\grqq. Dies ist mathematisch zwingend in den vorausgesetzten Konsistenzbedingungen der Arithmetik begründet: Eine Zulassung der Division durch Null würde die grundlegenden algebraischen Axiome verletzen und somit zu unvermeidbaren Widersprüchen innerhalb des Systems führen. Bei diesem Versuch hebt sich der Begriff also selbst auf und überschreitet die Grenze der logischen Möglichkeit \cite{kant_krv}.

Historisch sah bereits Bhāskara II. (12. Jh.) in der Division durch Null eine Annäherung an das Unendliche. In der modernen Analysis wird dieser Punkt als Singularität begriffen — ein Ort, an dem die üblichen Regeln der Logik ihre Gültigkeit verlieren \cite{fahse2014}. Die Entwicklung geht in mehrere Richtungen weiter. Einerseits bekommt die Null Bedeutungen, die nichts mit dem Nichts zu tun haben, wie der Referenzpunktaspekt. Andererseits werden die oben erwähnte algebraische Eigenschaften der Null in anderen Bereichen verallgemeinert, wie Nullvektor, Nullabbildung, Nullmatrix usw.

\section{Nichts und Die Leere Menge}

Die leere Menge ($\emptyset$ oder \{\})  formalisiert das Konzept der Absenz als ein eigenständiges mathematisches Objekt, dessen Existenz axiomatisch gesichert wird. Sie bildet das Fundament der Mengenlehre und kann somit als Grundbaustein der modernen Mathematik angesehen werden. Die leere Menge ist definiert \cite{kanamori2003} als die Menge, die keine Elemente enthält. Die entscheidende philosophische Frage \glqq Ist die leere Menge Nichts?\grqq{} muss aus mathematischer Sicht mit Nein beantwortet werden. Die leere Menge ist kein Nichts; sie ist ein Ding, das nichts enthält, oder eine Sammlung von keinen Objekten. 

Diese ontologische Differenzierung lässt sich philosophisch präzisieren, indem man aufzeigt, dass die leere Menge je nach Perspektive Facetten von zwei verschiedenen Kantischen Kategorien des Nichts vereint \cite{kant_krv}. Funktionell betrachtet lässt sie sich analog zu Kants \textit{ens imaginarium} (\glqq leere Anschauung ohne Gegenstand\grqq) begreifen: Ähnlich wie der reine Raum stellt sie das formale Gerüst (ein \glqq Gefä\ss\grqq) bereit, das in der Mathematik als existierendes Objekt behandelt wird, auch wenn es keine Elemente umschlie\ss{}t. 

Erkenntnistheoretisch und logisch betrachtet ist die leere Menge jedoch ein \textit{ens rationis} (\glqq leerer Begriff ohne Gegenstand\grqq). Sie ist ein widerspruchsfreies Gedankending, das zwar keine direkte Entsprechung in der empirischen Sinnenwelt besitzt, aber innerhalb des axiomatischen Systems als notwendiges Objekt vollumfänglich legitimiert ist. In beiden Lesarten ist sie in jedem Fall scharf von Kants \textit{nihil negativum} (dem logischen Unding) abzugrenzen, da ihre Konstruktion keinerlei Widerspruch enthält.

Sie ist in der Mengenlehre, insbesondere in Zermelo-Fraenkel (ZF), ein existierendes Objekt. Die Existenz der leeren Menge wird durch das Leermengenaxiom (oder durch Ableitung über das Separationsaxiom in Verbindung mit der Existenz einer anderen Menge) gesichert. Während die Existenz einer Universalmenge eher kontrovers ist, gibt es heute in der Mathematik so gut wie keine Kontroverse darüber, ob die leere Menge als Menge \cite{johnson1972} existiert.

Historisch erscheint die leere Menge recht spät auf der Bühne der Mathematik. Der erste, der über die leere Menge geschrieben hat, war George Boole \cite[Seite 28]{boole1854}, bei ihm hie\ss{} es noch Nichts. Er schreibt:

\begin{quote}
By a class is usually meant a collection of individuals, to each of which a particular name or description may be applied;
but in this work the meaning of the term will be extended so as to include the case in which but a single individual exists, answering to the required name or description, as well as the cases denoted by the terms \glqq nothing\grqq{} and \glqq universe,\grqq{} which as \glqq classes\grqq{} should be understood to comprise respectively \glqq no beings,\grqq{} \glqq all beings.\grqq
\end{quote}

Später betonte Peano (der das Symbol $\Lambda$ verwendete) durch seine Konstruktionen, dass es nur eine leere Menge gibt. Die Arbeiten von Peano sind sehr modern und faszinierend. Er hatte für die Notation $\Lambda$ einen logischen und ästhetischen Grund: Peano benutzte den Buchstaben $V$ (für Verum, das Wahre) als Symbol für die \glqq Allklasse\grqq{} (die Menge, die alles enthält) oder die logische Wahrheit. Da die leere Menge das logische Gegenteil der Allklasse ist (oder das \glqq Falsche\grqq{} in der Aussagenlogik), drehte er das Symbol einfach um\footnote{\cite[Seite XI]{Peano1889}. \glqq Signum $V$ significat verum, sive identitatem [...] Signum $\Lambda$ significat falsum, sive absurdum.\grqq}.

Aus $V$ wurde $\Lambda$. Er bezeichnete $\Lambda$ als die \glqq Nullklasse\grqq \footnote{\cite[Seite XI]{Peano1889}. \glqq Signum $\Lambda$ indicat classem, quae nullum continet individuum.\grqq}. Er formulierte die universelle algebraische Eigenschaften
$$A\cup \emptyset = A,\quad A\cap \emptyset = \emptyset,\quad  A\setminus A=\emptyset$$
für eine beliebige Menge $A$ (in moderner Notation) \cite[I.\S 6 Seite 12]{peano1894}, die eine wichtige Analogie zu den algebraischen Eigenschaften der Null aufweisen.

Die leere Menge wurde nicht allgemein akzeptiert, weil Mengen oft als Ansammlungen ihrer Elemente angesehen wurden. Dedekind versucht 1887 in seinem berühmten Werk \cite{dedekind2012} die Verwendung der leeren Menge (leeres System) bewusst zu vermeiden, er schreibt in $\S 1$: \glqq Dagegen wollen wir das leere System, welches gar kein Element enthält, aus gewissen Gründen hier ganz ausschlie\ss{}en, obwohl es für andere Untersuchungen bequem sein kann, ein solches zu erdichten\grqq. Die Arbeiten von Cantor und Hausdorff haben letztlich die leere Menge hoffähig gemacht. Es ist bis heute die grö\ss{}te didaktische Problematik mit der leeren Menge, dass durch die Fehlvorstellung wo Mengen als Ansammlungen angesehen werden, die leere Menge nicht als Menge akzeptiert \cite{fischbein1998} wird.

Es scheint aus der üblichen Definition der leeren Menge klar zu sein, dass es hier um eine Privation, eine Art echter Absenz geht: Die Menge aller Primzahlen zwischen 32 und 36 zeugt von der Absenz von Primzahlen mit dieser Eigenschaft. Hier offenbart sich jedoch eine tiefe philosophische Spannung zwischen der sprachlichen Vielfalt (Intension) und der mathematischen Eindeutigkeit (Extension). Wir können unendlich viele verschiedene \glqq negative Aussagen\grqq{} formulieren, die alle logisch wahr sind, aber inhaltlich gänzlich Verschiedenes beschreiben: Die \glqq Menge der runden Quadrate\grqq{} unterscheidet sich begrifflich klar von der \glqq Menge der rosafarbigen  Einhörner\grqq. In der Mengenlehre kollabieren diese semantischen Unterschiede jedoch durch das Extensionalitätsaxiom von Dedekind. Wie Quine \cite[Seite 3]{Quine2013} betont, werden Mengen nicht durch die Art ihrer Beschreibung (die \glqq offene Aussage\grqq), sondern ausschlie\ss{}lich durch ihre Elemente bestimmt. Da all diese negativen Beschreibungen keine Elemente liefern, referenzieren sie zwangsläufig auf exakt dasselbe Objekt. Wir sehen in dieser Reduktion keinen Mangel, sondern eine notwendige Klärung: Während Attribute oder Begriffsinhalte (Intensionen) oft vage Identitätskriterien haben, schafft die Mengenlehre durch die Reduktion auf die Extension klare ontologische Verhältnisse. Das Nichts in der leeren Menge verliert dadurch seinen relationalen Charakter und wird zu einem absoluten, einzigartigen Gegenstand ($\emptyset$), der alle inhaltlichen Unterschiede des Nichtseienden in sich absorbiert.

Die leere Menge hat spezifische logische Eigenschaften, die ihre aktive Rolle im formalen System belegen.
\begin{itemize}
\item Existenzaussagen über Elemente der $\emptyset$ sind stets falsch, da es kein Element gibt, das die Bedingung erfüllen könnte.
\item Allaussagen über Elemente der $\emptyset$ sind stets wahr (vacuously true). Zum Beispiel ist die Aussage \glqq Für alle Elemente $x$ der Menge $\emptyset$ gilt: $x$ ist blau\grqq{} wahr, da es keine Elemente gibt, die diese Bedingung verletzen könnten.
\end{itemize}

Das mathematische Nichts erzwingt somit eine logische Wahrheit, anstatt eine sinnlose Leerstelle zu sein. Ergänzend zur klassischen mengentheoretischen Perspektive lässt sich die funktionale Rolle der leeren Menge in der Sprache der Kategorientheorie präzisieren. In der Kategorie der Mengen ($\mathbf{Set}$) fungiert die leere Menge als das eindeutige {Initialobjekt} \cite{mac1998categories} (oder Anfangsobjekt). Diese universelle Eigenschaft besagt, dass für jede beliebige Menge $A$ genau eine einzige Abbildung $\emptyset \to A$ existiert — die sogenannte leere Funktion. Das mathematische Nichts ist somit kein isolierter Zustand, sondern zeichnet sich durch seine universelle Konnektivität zu allen Objekten des mathematischen Universums aus. Es bildet den absoluten strukturellen Ursprung, von dem aus jede Abbildung und damit jede Konstruktion möglich wird, ohne selbst auf inhaltliche Voraussetzungen angewiesen zu sein. Damit wird die leere Menge zum funktionalen \glqq Urpunkt\grqq, der in einer einzigartigen Beziehung zu jedem denkbaren \glqq Etwas\grqq{} steht.

Die leere Menge dient in der modernen Mathematik als absoluter Grundbaustein, aus dem sich weite Teile der formalen Disziplin konstruieren lassen. Dieser Aufbau beginnt bei den natürlichen Zahlen: Während die Peano-Axiome die Null als erstes Element\footnote{Das ist die moderne Formulierung, Peano hat noch die 1 gewählt, siehe \cite{Peano1889}, Seite 1., später aber 0, siehe \cite{peano1894}, Seite 28.} und eine Nachfolgerfunktion $S: x \rightarrow x+1$  noch axiomatisch als gegeben voraussetzen, wird dieses Fundament in der mengentheoretischen Konstruktion (z. B. nach John von Neumann) radikal aus der Leere erzeugt. Die Null wird hier direkt mit der leeren Menge identifiziert $0 = \emptyset$. Peanos vorausgesetzte Nachfolgerfunktion verschwindet als eigenständiges Axiom und wird durch die rein mengentheoretische Operation  $S(n):=n\cup \{n\}$ ersetzt \cite[Abschnitt 12]{Quine2013}. Durch diese Konstruktion werden die Peano-Axiome vollständig erfüllt, ohne dass Zahlen als ontologisch vorgegebene Objekte existieren müssen. Von diesem formalisierten Nichts aus – der leeren Menge – wird somit die gesamte Zahlenreihe generiert. Da sich wiederum nahezu die gesamte restliche Mathematik, einschlie\ss{}lich algebraischer Strukturen, Vektorräume und topologischer Räume, in die Mengenlehre einbetten lässt und Theoreme als Aussagen über Mengen betrachtet werden können, entspringt die Kette der mathematischen Existenz hier im strengen Sinne formal aus der Absenz.

\subsection{Das Nichts in der Informatik: Vom Systemabsturz zur strukturierten Absenz}

Die mathematische Unterscheidung zwischen der leeren Menge als eigenständigem Objekt und dem blo\ss{}en Fehlen von Information findet in der Informatik eine kritische und sehr praktische Entsprechung. Da Computer deterministische Maschinen sind, die auf wohldefinierten Speicherzuständen operieren, zwingt die Informatik dazu, das \glqq Nichts\grqq{} technisch exakt zu fassen. Dabei zeigen sich unterschiedliche Ausprägungen, deren Funktion und Nutzen sich stark unterscheiden.

Ein prominentes Beispiel für das unstrukturierte, destruktive Nichts ist die Problematik der \texttt{null}-Referenz, die ihr Erfinder Tony Hoare rückblickend als seinen \glqq Milliarden-Dollar-Fehler\grqq{} bezeichnete\footnote{https://www.infoq.com/presentations/Null-References-The-Billion-Dollar-Mistake-Tony-Hoare/ (abgerufen am 23.1.2026)}. In vielen imperativen Programmiersprachen wird \texttt{null} (oder \texttt{nil}) als Zeiger auf einen ungültigen Speicherbereich verwendet. Dieses Nichts besitzt keine Eigenschaften und keine Struktur. Versucht das Programm, auf eine Eigenschaft dieses Nichts zuzugreifen (Dereferenzierung), führt dies -- analog zur Division durch Null in der Arithmetik -- in der Regel zu einem sofortigen Systemabsturz. Das System kollabiert, weil es versucht, auf Struktur zuzugreifen, wo absolute Leere herrscht.

Ein funktional völlig anderes Konzept ist die leere Datenstruktur, wie etwa ein leeres Array \texttt{[]} oder eine leere Zeichenkette (\texttt{""}). Im Gegensatz zum \texttt{null}-Zeiger ist ein leeres Array ein existierendes Objekt im Speicher. Es besitzt eine Typ-Signatur und Methoden, aber seine Länge ist null. Es ist die direkte informatische Umsetzung der leeren Menge ($\emptyset$): Ein Programm kann fehlerfrei über ein leeres Array iterieren -- die Schleife wird schlicht null Mal ausgeführt. Das Nichts ist hier in ein formales Gefä\ss{} gegossen, mit dem das System sicher operieren kann.

Eine weitere Ausprägung findet sich in relationalen Datenbanken (SQL) \cite{codd1979}. Hier fungiert \texttt{NULL} nicht zwingend als ontologische Leere, sondern oft als epistemisches Nichts: Es steht für \glqq unbekannt\grqq{} oder \glqq nicht zutreffend\grqq. Diese funktionale Bedeutung zwingt Datenbanksysteme zur Verwendung einer dreiwertigen Logik (wahr, falsch, unbekannt). Der Vergleich \texttt{NULL = NULL} ergibt in SQL nicht \textit{wahr}, sondern wiederum \textit{unbekannt}, da zwei unbekannte Werte nicht logisch zwingend identisch sind. Das Nichts wird hier genutzt, um fehlendes Wissen systemisch abzubilden.

Um die Fehleranfälligkeit von \texttt{null}-Referenzen zu beheben, nutzen moderne funktionale Programmiersprachen das Konzept der \texttt{Option}- oder \texttt{Maybe}-Typen \cite{pierce2002}. Hierbei wird das Nichts (\textit{None} oder \textit{Nothing}) als ein wohlstrukturierter Zustand innerhalb eines Datentyps begriffen. Der Compiler zwingt den Programmierer, den Fall der Abwesenheit eines Wertes explizit zu behandeln. Das Programm stürzt nicht ab, sondern verarbeitet die Information über die Abwesenheit eines Wertes als validen logischen Pfad. Das Nichts verliert seine destruktive Kraft und wird sicher in den Kontrollfluss eingebettet.

In der tieferliegenden Typentheorie und formalen Semantik wird die Unterscheidung zwischen verschiedenen Formen des Nichts noch präziser gefasst. Hier stehen Typen wie der \texttt{Unit}-Typ (in C oder Java oft als \texttt{void} bezeichnet) Berechnungen gegenüber, die den \texttt{Bottom-Typ} ($\bot$) aufweisen \cite{pierce2002}. Eine Funktion mit dem Rückgabetyp \texttt{Unit} führt eine Aktion aus, liefert aber keinen informativen Wert zurück -- der Typ enthält exakt einen einzigen Wert (ein Informationsgehalt von $0$ Bit). Der \texttt{Bottom-Typ} hingegen modelliert Berechnungen, die niemals ein Ergebnis liefern, beispielsweise weil sie in einer Endlosschleife gefangen sind oder mit einem Fehler abbrechen. Er ist der Typ, der absolut keine Werte enthält, und ist somit das direkteste formale Analogon zur echten leeren Menge.

Diese Differenzierungen in der Informatik unterstreichen die Hauptthese dieses Kapitels: Wo das Nichts nicht als strukturiertes Objekt ($\emptyset$), sondern als blo\ss{}e, ungültige Leerstelle begriffen wird, verliert das System seine Konsistenz. Sobald die Absenz jedoch als typisierte Struktur begriffen wird, lässt sich mit ihr sicher und produktiv operieren.

\section{Nichts ist nützlich}\label{sec:nutz}

Bisher haben wir das Nichts als ontologischen Status (Existiert es? Wie?) oder als axiomatisches Fundament (Mengenlehre) betrachtet. Doch in der mathematischen Praxis ist das Nichts oft weit mehr als eine passive Leerstelle: Es ist ein aktives, operatives Werkzeug. Das Nichts wird \glqq benutzt\grqq, um Strukturen überhaupt erst sichtbar zu machen. Wir möchten das hier an drei Beispielen illustrieren.

Brian Rotman zieht in seiner semiotischen Untersuchung \cite{rotman1987} eine faszinierende Parallele zwischen der Einführung der mathematischen Null und der Erfindung des Fluchtpunkts in der perspektivischen Malerei der Renaissance. Der Fluchtpunkt ist jener Punkt im Bild, an dem sich alle parallelen Linien treffen und in der Unendlichkeit verschwinden. Physisch ist er auf der Leinwand ein \glqq Nichts\grqq, ein Punkt ohne Ausdehnung. Doch strukturell ist er das organisierende Zentrum, das dem Bild erst seine Tiefe und Räumlichkeit verleiht. Ohne dieses visualisierte Nichts bricht die Illusion des dreidimensionalen Raumes zusammen. Rotman argumentiert, dass die Null in der Mathematik eine ähnliche Funktion erfüllt: Sie ist ein \glqq Meta-Zeichen\grqq, das den Code der Zahlen organisiert, so wie der Fluchtpunkt die visuelle Szene organisiert.

Diese organisierende Kraft zeigt sich nirgends deutlicher als in der Algebra, wo das Nichts als kreativer Trick eingesetzt wird. Ein Standardverfahren in der Analysis und Algebra besteht darin, einen Term künstlich zu komplizieren, indem man eine \glqq nützliche Null\grqq{} addiert. Formal lässt sich dies als $a = a + b - b$ ausdrücken. Man fügt etwas hinzu und nimmt es sofort wieder weg. Der Wert des Ausdrucks bleibt unverändert, doch seine Gestalt wandelt sich, wodurch verborgene Strukturen erkennbar werden. Ein klassisches Beispiel hierfür ist die quadratische Ergänzung. Um eine Gleichung der Form $x^2 + px$ in eine binomische Formel zu verwandeln, addiert man das Quadrat der halben Vorzahl und subtrahiert es sofort wieder:

$$x^2 + px = \underbrace{x^2 + px + \left(\frac{p}{2}\right)^2}_{\text{Struktur}} - \left(\frac{p}{2}\right)^2 = \left(x + \frac{p}{2}\right)^2 - \left(\frac{p}{2}\right)^2$$

Durch das Hinzufügen dieses \glqq komplizierten Nichts\grqq{}

$$0 = +\left(\frac{p}{2}\right)^2 - \left(\frac{p}{2}\right)^2 $$

wird die Lösungsstruktur der Gleichung plötzlich sichtbar. Das Nichts fungiert hier nicht als Mangel, sondern als Katalysator für Erkenntnis. Dieses Prinzip der produktiven Null lässt sich auch in der Physik, insbesondere bei den Erhaltungssätzen, wiederfinden. Die Kirchhoffschen Regeln der Elektrotechnik fordern etwa, dass die Summe aller Ströme in einem Knotenpunkt Null ergibt ($\sum I=0$). Die Null markiert hier keine Leere, sondern ein dynamisches Gleichgewicht. Diese Logik der Balance zeigt sich ebenso in den Nullsummenspielen der Spieltheorie. Betrachtet man zudem das Noether-Theorem, das Erhaltungssätze mit Invarianzen verknüpft, so wird die Null zum formalen Ausdruck fundamentaler Symmetrien eines Systems. In diesem Sinne ist das mathematische Nichts oft weniger ein \glqq Loch\grqq{} im Sein, sondern vielmehr der operative Dreh- und Angelpunkt, um den sich die mathematischen Strukturen bewegen. Es ist, wie Rotman andeutet, der Punkt, von dem aus das System sich selbst organisiert.

Die mathematische Intuition der Null als Zustand perfekten Gleichgewichts findet in der modernen Quantenfeldtheorie eine spektakuläre Bestätigung. Hier wird das physikalische Vakuum — das absolute Nichts der klassischen Mechanik — als ein dynamisches Feld begriffen. Das Quantenvakuum ist nicht leer, sondern wird als der Zustand niedrigster Energie ($|0\rangle$) definiert, der jedoch beständigen Fluktuationen unterliegt \cite{visser1995}. Diese sogenannten Vakuumfluktuationen führen zur spontanen Entstehung virtueller Teilchen-Antiteilchen-Paare, die sich unmittelbar wieder vernichten. Dieser physikalische Sachverhalt korrespondiert auf faszinierende Weise mit Platons relationalem Nichts aus dem {Sophistes}: Das Vakuum ist nicht das blo\ss{}e Gegenteil des Seins, sondern eine \glqq Andersheit\grqq, die als unerschöpfliches Reservoir an Möglichkeiten fungiert \cite{visser1995}. Das Nichts ist hier im wahrsten Sinne des Wortes \glqq schwanger\grqq{} mit Struktur. Operativ wird dieses \glqq aktive Nichts\grqq{} etwa im Casimir-Effekt messbar, bei dem die blo\ss{}e Anwesenheit von Vakuumfeldern zwischen zwei Metallplatten eine messbare Kraft ausübt. Damit erweist sich das Nichts auch in der Physik nicht als Mangel, sondern als ein hochenergetischer Dreh- und Angelpunkt der Realität.

\section{Nichts und die Essenz der Mathematik}

\subsection{Der sokratische Nullpunkt: Erkenntnis aus dem Nichtwissen}

Der Übergang von der formalen Struktur der Mathematik zu ihrer \glqq Essenz\grqq{} im menschlichen Geist erfordert eine Untersuchung des Zustands, der jeder Entdeckung vorausgeht: das bewusste Nichtwissen. In der Geschichte der Philosophie fungiert Sokrates als die personifizierte \glqq Nullstelle\grqq{} des Denkens. Sein berühmtes Bekenntnis zum Nichtwissen ist nicht als resignative Geste zu verstehen, sondern als die Etablierung eines erkenntnistheoretischen Vakuums, das erst den Raum für die Entfaltung logischer Notwendigkeit schafft.

Der Ausgangspunkt des sokratischen Dialogs ist in der Regel das unhinterfragte Scheinwissen der Gesprächspartner. Die Aporie (gr. aporía – Weglosigkeit) ist folglich nicht der Beginn, sondern vielmehr das erste gro\ss{}e Ziel und der entscheidende Durchgangspunkt der sokratischen Untersuchung. Wenn das anfängliche Scheinwissen kollabiert, entsteht ein Zustand der Ratlosigkeit, ein geistiges \glq Nichts\grq{}. Genau diesen Zustand beschreibt Karl Jaspers \cite{jaspers1997} als \glq Grenzsituation des Geistes\grq, in der Sokrates selbst als \glq leere Mitte\grq{} fungiert. Er besitzt kein fertiges Lehrgebäude (kein \glqq Etwas\grqq), sondern schafft durch sein Fragen einen Raum der Maieutik (Hebammenkunst). Die Aporie ist somit kein blo\ss{}er Mangel, sondern jener produktive sokratische Nullpunkt, an dem das Denken über sich selbst hinausweist und von dem aus die Konstruktion der Wahrheit ex nihilo aus der reinen Vernunft beginnen kann.

Im Menon-Dialog \cite{ebert2018} (82b–85b) wird dies am Beispiel des Sklaven explizit: Bevor dieser die geometrische Wahrheit über die Verdoppelung eines Quadrats erkennen kann, muss sein anfängliches Scheinwissen dekonstruiert werden. Er wird durch gezieltes Fragen in den Zustand der Ratlosigkeit geführt – einen Zustand des geistigen \glqq Nichts\grqq{} in Bezug auf das Problem, den Sokrates mit der heilsamen Betäubung durch einen Zitterrochen vergleicht.

Um die schöpferische Kraft dieses \glqq geistigen Nichts\grqq{} zu verstehen, ist eine genauere Analyse des mathematischen Objekts im {Menon} aufschlussreich. Der Sklave beginnt den Dialog mit einer plausiblen, aber falschen Konstruktion: Er vermutet, dass die Verdoppelung der Quadratfläche durch eine Verdoppelung der Seitenlänge erreicht wird: ein Irrtum, der zu einer vervierfachten Fläche führt. Die Aporie tritt exakt in dem Moment ein, in dem dieses Scheinwissen kollabiert. Dieses \glqq Nichts\grqq{} der Gewissheit ist jedoch hochgradig strukturiert: Es ist die Grenze, an der die geometrische Wahrheit als diagonale Konstruktion erst sichtbar werden kann.

Diesen Gedanken des produktiven Nichtwissens und der schöpferischen Negation greifen moderne Mathematiker und Philosophen auf, um die Natur mathematischer Entdeckung zu beschreiben. Alfréd Rényi nutzt in seinen Dialogen \cite{renyi1967} die sokratische Form, um zu demonstrieren, dass mathematische Objekte keine greifbaren Realitäten sind, sondern aus dem Prozess des Hinterfragens entstehen. In Rényis fiktivem Dialog mit Sokrates wird deutlich, dass die Mathematik ihre Objekte aus der Leere der Abstraktion gewinnt; sie benötigt keine materielle Basis, sondern lediglich die Reinheit des logischen Raums, der durch das Eingeständnis des Nichtwissens geöffnet wurde.

Einen weiteren Schritt geht Imre Lakatos \cite{lakatos2013}, indem er die sokratische Methode direkt auf die mathematische Forschungspraxis überträgt. In seiner modernen Lesart betont er, dass die mathematische Entdeckung kein linearer Aufbau ist, sondern ein dialektischer Sprung, der die Widerlegung als Sprungbrett nutzt. Für Lakatos wächst mathematisches Wissen nicht durch das Anhäufen unumstö\ss{}licher Fakten, sondern durch einen dialektischen Prozess von Vermutungen und Negationen (Kritik). Das \glqq Nichts\grqq{} einer widerlegten Theorie oder eines Gegenbeispiels ist hierbei der entscheidende Motor: Erst durch das Aufzeigen dessen, was nicht der Fall ist, schärft sich die mathematische Struktur. Lakatos zeigt, dass die Mathematik eine \glqq Logik der Entdeckung\grqq{} besitzt, die auf der ständigen Überwindung des Irrtums beruht. Die mathematische Wahrheit \glqq erwacht\grqq{} nicht trotz, sondern wegen des Vakuums, das durch die Zerstörung falscher Intuitionen entstanden ist; sie wird somit nicht passiv gefunden, sondern in einem aktiven Prozess aus dem \glqq Nichts\grqq{} der Widerlegung heraus gemei\ss{}elt.

Die Analyse dieses Prozesses zeigt eine faszinierende Parallele zur von Neumannschen Konstruktion der Zahlen. So wie die Mengenlehre die gesamte Arithmetik aus der leeren Menge $\emptyset$ entfaltet, so entfaltet der sokratische Prozess komplexe Wahrheiten aus einem Geist, der scheinbar leer von Fachwissen ist. Mathematisch gedeutet bedeutet dies: Die logischen Strukturen sind bereits als Potenzialität im \glqq Nichts\grqq{} des Geistes vorhanden. Das Nichtwissen ist hierbei kein Mangel, sondern die aktive Null, die das System in Bewegung setzt. In dieser Lesart ist die Leere des Geistes die notwendige Voraussetzung für die creatio ex nihilo mathematischer Gewissheit, wie sie in den nachfolgenden Abschnitten bei Nikolaus von Kues und János Bólyai zur zentralen Kategorie wird.

\subsection{Das belehrte Nichtwissen als schöpferisches Prinzip}

Nikolaus von Kues (1401–1464) markiert mit seinem Werk einen fundamentalen Perspektivwechsel in der Mathematikphilosophie, indem er die Gegenstände der Mathematik als freie Schöpfungen des menschlichen Geistes (mens humana) aus dem Nichts auffasst. Dieser Ansatz steht in direktem Gegensatz sowohl zur platonischen als auch zur aristotelischen Tradition. In der antiken Philosophie galt die Mathematik einhellig als rein passive Schau oder Abstraktion. Aristoteles fasst diese reduktive Sichtweise in seiner Metaphysik prägnant zusammen: 

\begin{quote}
\glqq Der Mathematiker stellt nun Betrachtungen an über das, was aus einer Wegnahme (Aphairesis) hervorgeht. Er betrachtet nämlich die Dinge, indem er alles Sinnliche weglässt [...] er lässt nur das Quantum und das Zusammenhängende übrig.\grqq{} (Met. K 3)
\end{quote}

Bei Cusanus hingegen wird das Erfassen mathematischer Gegenstände zu einem aktiven geistigen Konstruieren (praxis) anstelle einer passiven theoria, siehe Nickel \cite{nickel2021}. Wie David Albertson \cite{albertson2014} und Fritz Nagel \cite{nagel2007} aufzeigen, formt Cusanus die Mathematik damit zu einer \glqq theologia mathematica\grqq, in der das Mathematisieren als schöpferische Wissenschaft betrieben wird.

Cusanus etabliert in diesem Rahmen eine tiefgreifende Analogie zwischen der göttlichen und der menschlichen Schöpfung. Diese theologische Metapher beginnt beim Punkt, der dem Nichts am nächsten ist. In seinem Werk \textit{De theologicis complementis} bedient er sich einer rein metaphorisch zu verstehenden \glqq Schöpfungsgeschichte\grqq{} (aenigma), in der endliche geometrische Konstruktionen als Symbole fungieren, um die metaphysische creatio ex nihilo zu spiegeln. In der Übersetzung von Gregor Nickel \cite{nickel2007} formuliert Cusanus wörtlich:

\begin{quote}
  \glqq Der Schöpfer scheint also zwei Dinge geschaffen zu haben, nämlich [einerseits] nahe beim Nichts den Punkt — zwischen Punkt und Nichts gibt es nämlich kein Mittleres [...] — und andererseits nahe bei sich das Eine. Und diese fügte er zusammen, so da\ss{} [es] Ein Punkt sei. In diesem einen Punkt war die Einfaltung des Universums.\grqq{}
\end{quote}

Analog zu dieser metaphysischen Entfaltung (Explicatio) beschreibt Cusanus, wie der menschliche Geist die mathematischen Gegenstände konstruiert. Wenn der Geist eine Figur darstellen will, beginnt er laut Cusanus mit einem einzigen Punkt, erweitert diesen zur Linie und biegt ihn zu Winkeln ab, um die Fläche einzuschlie\ss{}en und das Vieleck zu vollenden. Auf diese Weise imitiert der menschliche Geist den göttlichen Schöpfungsakt im eigenen Denken. Dass die Mathematik hierbei nicht der empirischen Natur entnommen wird (wie noch bei der aristotelischen Aphairesis), sondern eine primäre Schöpfung des Geistes ist, belegt Cusanus in De beryllo unmissverständlich:

\begin{quote}
 \glqq [D]ass unser Geist, der die mathematischen Gegenstände schafft, das was er schaffen kann, wahrer und wirklicher in sich hat, als es au\ss{}er ihm ist. (...) Und so ist es bei allem dergleichen, beim Kreis, bei der Linie, beim Dreieck, auch bei unserem Zahlbegriff, kurz bei allem, was seinen Ursprung aus dem menschlichen Geist nimmt und der Natur entbehrt.\grqq{} (De beryllo c. 56)
\end{quote}

Die Grundlage für diese schöpferische Kraft liegt in der Gottebenbildlichkeit (imago dei) des Menschen. Für Cusanus steht daher fundamental fest, dass \glqq nichts vornehmer als der menschliche Geist\grqq{} sei. In seiner Freiheit setzt der Geist den Dingen Grenzen und bestimmt mathematische Gegenstände aus sich heraus \cite{nickel2007}:

\begin{quote}
  \glqq Der menschliche Geist, der ein Bild des absoluten Geistes ist, setzt in seiner menschlichen Freiheit allen Dingen in seinem Denken Grenzen, weil der Geist mit seinen Begriffen alles ausmi\ss{}t. Er setzt eine Grenze für die Linien, macht sie lang oder kurz, und setzt so viele Begrenzungspunkte in ihnen, wie er will.\grqq{} (De venatione sapientiae c. 27: h XII, N. 82 Z. 13-17)
\end{quote}

Gerade weil der menschliche Geist den mathematischen Gegenstand gänzlich aus sich heraus entwickelt, kann ihm sein Gegenstand in grö\ss{}ter Genauigkeit bekannt sein. Die daraus resultierende grö\ss{}tmögliche Sicherheit der mathematischen Erkenntnis ist somit das Resultat einer aktiven Leistung des menschlichen Verstandes (ratio) und markiert den endgültigen Bruch mit der Antike \cite{nickel2018}. Zusammengefasst: So wie Gott die Welt aus dem absoluten Nichts erschaffen hat, erschafft der menschliche Geist neue Mathematik aus dem relationalen Nichts\footnote{vgl. zur weiteren Vertiefung des cusanischen Mathematik- und Unendlichkeitsverständnisses auch Müller \cite{muller2010} sowie Pukelsheim/Schwaetzer \cite{pukelsheim2005}.}.

\subsection{Die Freiheit der Konstruktion: Von Bólyai zu Cantor}

Die von Nikolaus von Kues theoretisch begründete schöpferische Freiheit des Geistes fand im 19. Jahrhundert ihre radikale mathematische Realisierung. Dieser Abschnitt untersucht den Moment, in dem die Mathematik die Fesseln der physikalischen Anschauung abstreift und das \glqq Nichts\grqq{} als Raum unbegrenzter Möglichkeiten besetzt.

Der wohl berühmteste Ausspruch über die mathematische Schöpfungskraft stammt von János Bólyai, der die folgende berühmte Textpassage in einem Brief an seinen Vater, Wolfgang Bólyai, am 3. November 1823 über seine Entdeckung \cite[Seite 187]{barna1970} der nichteuklidischen Geometrie geschrieben\footnote{Übersetzung von \cite{stwolfgang}, Seite 85.} hat:

\begin{quote}
	wenn Sie, mein teurer Vater, es sehen werden, so werden Sie es erkennen;
jetzt kann ich nichts weiter sagen, nur so viel: da\ss{} ich \emph{aus nichts eine neue, andere Welt geschaffen habe.} Alles, was ich bisher geschickt habe, ist ein Kartenhaus im Vergleich zu einem Turme.
\end{quote}

Mit der Entdeckung einer nicht-euklidischen Geometrie bewies Bólyai, dass mathematische Systeme nicht durch die physische Welt (den euklidischen Raum) vorgegeben sind. Indem er das Parallelenaxiom nicht als \glqq Wahrheit\grqq, sondern als wählbare Setzung begriff, wurde das Nichts – die Abwesenheit einer zwingenden äu\ss{}eren Vorgabe – zum Fundament einer neuen Geometrie. Wie Jeremy Gray \cite{gray2007} betont, markiert dieser Moment den Übergang von der Mathematik als \glqq Abbild der Natur\grqq{} zur Mathematik als \glqq freies Gedankenkonstrukt\grqq.

Dieser Weg mündet konsequent in das Werk von Georg Cantor, der feststellte \cite[Seite 182]{cantor1932}: 

\begin{quote}
	... das Wesen der Mathematik liegt gerade in ihrer Freiheit.
\end{quote}

Für Cantor ist die Mathematik nicht mehr an das Sein materieller Objekte gebunden, sondern nur noch an die logische Widerspruchsfreiheit ihrer Konstruktionen. Mit der Erfindung der Mengenlehre gab Cantor dem Nichts eine explizite Form: Die leere Menge $\emptyset$ ist nicht nur ein Symbol für Abwesenheit, sondern der aktive Ausgangspunkt, aus dem durch rein geistige Akte (Bündelung, Paarbildung) die Unendlichkeit der transfiniten Zahlen generiert wird. Die Mathematik wird hier zur absoluten Schöpfung ex nihilo, die keine Rechtfertigung au\ss{}erhalb ihrer eigenen Konsistenz benötigt.

Eine radikale formallogische Umsetzung dieser Schöpfung aus dem Nichts liefert im 20. Jahrhundert George Spencer Brown in seinem Werk Laws of Form. Wie Martin Rathgeb \cite{rathgeb2016} in seiner detaillierten Analyse aufzeigt, konstruiert Spencer Brown eine fundamentale \glqq Nussschalenmathematik\grqq{} aus einem absoluten Nullpunkt. Der Ausgangspunkt seines sogenannten Indikationenkalküls ist die absolute Leere, die Brown später als das \glqq Ock\grqq{} bezeichnete: ein Raum, der derart abwesend ist, \glqq da\ss{} es überhaupt nichts ist, nicht einmal leer.\grqq{} Die Erschaffung mathematischer Strukturen erfolgt hier durch einen einzigen schöpferischen Primärakt: die Unterscheidung (Distinktion). Das Ziehen einer Grenze generiert in dieser Leere zwei Zustände (einen markierten und einen unmarkierten) und verleiht dem Nichts allererst eine Form. Aus diesen elementaren Konzepten der Unterscheidung (distinction) und Bezeichnung (indication) generiert Spencer Brown in der Folge systematisch seine \glqq Primäre Arithmetik\grqq{} und darauf aufbauend eine \glqq Primäre Algebra.\grqq{} Sein Kalkül liefert damit den modernen logischen Beleg dafür, wie eine gesamte Mathematik ex nihilo aus der reinen Setzung einer Unterscheidung entfaltet werden kann.

Den zeitgenössischen Abschluss dieser Entwicklung bildet das Konzept der Mathematik als Wissenschaft von Mustern, wie es Erich Ch. Wittmann \cite{wittmann2018} präzisiert. Wittmann definiert Mathematik nicht als eine Ansammlung isolierter Sätze, sondern als die systematische Untersuchung von Strukturen und Regelmä\ss{}igkeiten, die er als \glqq Muster\grqq{} bezeichnet. Diese Sichtweise radikalisiert den Gedanken der schöpferischen Freiheit: Mathematische Muster sind für Wittmann keine passiven Entdeckungen in einer materiellen Welt, sondern Konstrukte, die durch \glqq operative Handlungen\grqq{} des Geistes entstehen. Hier schlie\ss{}t sich der Kreis zum Thema des Nichts: Die Bausteine der Mathematik haben keine substanzielle Materie; sie sind rein relationale Gefüge. In der Wissenschaft von Mustern wird deutlich, dass die Mathematik ihre Gegenstände nicht vorfindet, sondern sie durch die Setzung von Beziehungen und Regeln erst erschafft. Wittmann betont dabei die Bedeutung des \glqq aktiven Lernens\grqq{} und des \glqq eigenständigen Konstruierens\grqq. Für uns bedeutet dies: Die Freiheit, von der Cantor sprach, wird bei Wittmann zur didaktischen und fachlichen Notwendigkeit. Wenn die Mathematik die Wissenschaft von Mustern ist, dann ist sie das Studium jener Ordnungen, die der menschliche Geist in den leeren Raum des Denkbaren hineinwebt. Das Nichts ist somit nicht das Ende der Mathematik, sondern der absolut notwendige Freiheitsgrad, ohne den die Konstruktion neuer Muster nicht möglich wäre.

Diese mathematische Freiheit, Muster in einem substanzlosen Raum zu entwerfen, findet ihre philosophische Erdung in Kants Begriff des \textit{ens imaginarium}. Für Kant ist der reine Raum eine \glqq leere Anschauung ohne Gegenstand\grqq{}: eine blo\ss{} formale Bedingung, die zwar selbst kein Objekt ist, aber erst den Rahmen für jede mögliche Konstruktion bereitstellt \cite{kant_krv}. 

In der modernen Lesart Wittmanns wird dieses kantische \glqq Einbildungsding\grqq{} zum aktiven Spielfeld der Vernunft: Das Nichts ist nicht die Abwesenheit von Ordnung, sondern die reine, formale Potenzialität, in der der menschliche Geist seine Strukturen \textit{ex nihilo} entfaltet. Damit schlie\ss{}t sich der Kreis von Kants transzendentaler Ästhetik zur absoluten Freiheit der modernen Mathematik: Das Nichts bildet die notwendige \glqq leere Form\grqq{}, ohne die der schöpferische Akt der Musterbildung keinen Raum zur Entfaltung fände.

\section{Fazit: Das Nichts als Raum schöpferischer Freiheit}

Das Nichts erweist sich in der Philosophie der Mathematik nicht als Abgrund der Bedeutungslosigkeit, sondern als jener produktive Raum, in dem sich Abwesenheit durch den Akt der Freiheit in formale Existenz verwandelt. Die Untersuchung zeigt, dass die Mathematik primär mit dem relationalen Nichts operiert – einem Prinzip der Andersheit, das die notwendige Bedingung für Pluralität und Struktur darstellt. Diese schöpferische Kraft manifestiert sich in den zwei zentralen Säulen der mathematischen Formalisierung:

\begin{itemize}
	\item \textbf{Die Zahl Null (0):} Ihre Entwicklung vom blo\ss{}en Platzhalter hin zum aktiven Fundament markiert den ersten gro\ss{}en Befreiungsschlag des mathematischen Denkens. Als neurobiologische Kategorie erfordert sie ein aktives neuronales Feuern, das die Leere als quantitative Grö\ss{}e begreift. Sie ist nicht nur das neutrale Element der Arithmetik, sondern fungiert als ontologische Symmetrieachse, die das System erst in seiner Gesamtheit als Gleichgewicht von Sein und Gegen-Sein ordnet.
	\item \textbf{Die leere Menge ($\emptyset$):} Sie ist das axiomatische \glqq Etwas\grqq, das als Ur-Objekt die gesamte Kette der mathematischen Existenz generiert. In ihrer Rolle als kategorientheoretisches Initialobjekt erweist sie sich als funktionaler Startpunkt jeder Abbildung. Dass diese formale Strukturierung des Nichts zwingend ist, zeigt sich ex negativo in der Informatik: Wo das Nichts als unstrukturierte Leerstelle (Null-Referenz) verbleibt, kollabiert das System; erst die Begreifung der Absenz als strukturierter Informationstyp (Option-Typ) wahrt die logische Konsistenz.
\end{itemize}

Diese Idee der Konstruktion zieht sich von den formalen Objekten bis hin zur Erkenntnistheorie. Die Analyse der sokratischen Aporie verdeutlicht, dass mathematische Wahrheit nicht passiv vorgefunden wird. Das bewusste Nichtwissen schafft erst jenes schöpferische Vakuum, in dem durch die Destruktion falscher Vorannahmen die Vernunft ihre eigene Struktur \textit{ex nihilo} entfalten kann. In der Tradition von Nikolaus von Kues vollendet sich dieser Gedanke: Der menschliche Geist erschafft Mathematik als freie Schöpfung aus dem relationalen Nichts. Mit Bólyai und Cantor emanzipiert sich diese Freiheit endgültig von der physischen Anschauung, indem sie das Nichts zum souveränen Spielfeld der Konstruktion erklärt. George Spencer Brown radikalisiert diesen Ansatz im 20. Jahrhundert formallogisch, indem er demonstriert, wie sich aus der absoluten Leere allein durch den primären Akt der Unterscheidung ein vollständiger Kalkül entfalten lässt. In der zeitgenössischen Definition der Mathematik als \glqq Wissenschaft von Mustern\grqq{} findet diese Entwicklung schlie\ss{}lich ihre Erfüllung: Mathematik ist das Studium jener Beziehungsgefüge, die der menschliche Geist in den leeren Raum des Denkbaren hineinwebt. Das Nichts ist somit nicht das Ende der Mathematik, sondern die unentbehrliche Bedingung ihrer Freiheit und der Ursprung jeder neuen Welt.

\printbibliography

@book{metzler,
	title={Metzler Lexikon Philosophie: Begriffe und Definitionen},
	author={Prechtl, P. and Burkard, F.P.},
	isbn={9783476054692},
	year={2015},
	publisher={J.B. Metzler},
	address={Stuttgart},
	doi={10.1007/978-3-476-05469-2}
}

@book{albertson2014,
  title={Mathematical theologies: Nicholas of Cusa and the legacy of Thierry of Chartres},
  author={Albertson, David},
  year={2014},
  publisher={Oxford University Press}
}

@article{nagel2007,
  title={{Nicolaus Cusanus -- mathematicus theologus. Unendlichkeit und Infinitesimalmathematik}},
  author={Nagel, F},
  journal={Trierer Cusanus Lecture},
  volume={13},
  year={2007}
}

@article{fahse2014,
	title={Vorstellungen zur {Null} im {Kontext} der {Division} durch {Null}},
	author={Fahse, Christian},
	journal={mathematica didactica},
	volume={37},
	number={1},
	pages={5--29},
	year={2014},
	doi={10.18716/ojs/md/2014.1118}
}

@book{neugebauer1969,
	title={The exact sciences in antiquity},
	author={Neugebauer, Otto},
	volume={9},
	year={1969},
	publisher={Courier Corporation},
	address={New York}
}

@article{neugebauer1941,
	title={On a special use of the sign \glqq zero\grqq{} in cuneiform astronomical texts},
	author={Neugebauer, Otto},
	journal={Journal of the American Oriental Society},
	volume={61},
	number={4},
	pages={213--215},
	year={1941},
	doi={10.2307/594251}
}

@phdthesis{chrisomalis2003,
	title={The comparative history of numerical notation},
	author={Chrisomalis, Stephen},
	year={2003},
	school={McGill University},
	address={Montreal}
}

@book{ebert2018,
  author    = {Platon},
  title     = {Menon. Übersetzung und Kommentar von Theodor Ebert},
  volume    = {134},
  publisher = {Walter de Gruyter GmbH \& Co KG},
  year      = {2018},
  address   = {Berlin},
  doi       = {10.1515/9783110502121}
}

@book{jaspers1997,
	title={Die {ma\ss gebenden} Menschen: Sokrates, Buddha, Konfuzius, Jesus},
	author={Jaspers, Karl},
	year={1997},
	publisher={Piper},
	address={M{\"u}nchen}
}

@book{kaplan2000,
	title={The nothing that is: A natural history of zero},
	author={Kaplan, Robert},
	year={1999},
	publisher={Oxford University Press},
	address={Oxford}
}

@article{nieder2016,
	author  = {Nieder, Andreas},
	title   = {Representing Something Out of Nothing: The Dawning of Zero},
	journal = {Trends in Cognitive Sciences},
	year    = {2016},
	volume  = {20},
	number  = {11},
	pages   = {830--842},
	doi     = {10.1016/j.tics.2016.08.008}
}

@article{kutter2024,
	title={Single-neuron representation of nonsymbolic and symbolic number zero in the human medial temporal lobe},
	author={Kutter, Esther F and Dehnen, Gert and Borger, Valeri and Surges, Rainer and Nieder, Andreas and Mormann, Florian},
	journal={Current Biology},
	volume={34},
	number={20},
	pages={4794--4802},
	year={2024},
	publisher={Elsevier},
	doi={10.1016/j.cub.2024.08.041}
}

@book{mac1998categories,
	title={Categories for the working mathematician},
	author={Mac Lane, Saunders},
	volume={5},
	series={Graduate Texts in Mathematics},
	year={1998},
	publisher={Springer Science \& Business Media},
	address={New York},
	doi={10.1007/978-1-4757-4721-8}
}

@book{renyi1967,
	author = {Alfréd Rényi},
	publisher = {Holden-Day},
	title = {Dialogues on Mathematics},
	year = {1967},
	address = {San Francisco}
}

@book{lakatos2013,
	title={Beweise und Widerlegungen: die Logik mathematischer Entdeckungen},
	author={Lakatos, Imre},
	year={2013},
	publisher={Springer-Verlag},
	address={Berlin}
}

@book{pierce2002,
	author    = {Pierce, Benjamin C.},
	title     = {Types and Programming Languages},
	publisher = {MIT Press},
	year      = {2002},
	isbn      = {978-0262162098},
	address   = {Cambridge, MA}
}

@book{visser1995,
	author    = {Visser, Matt},
	title     = {Lorentzian Wormholes: From Einstein to Hawking},
	publisher = {American Institute of Physics},
	year      = {1995},
	address   = {Woodbury, NY},
	note      = {Besonders Kapitel 4 über das Quantenvakuum und Fluktuationen.}
}

@book{gray2007,
	title={Worlds out of Nothing: A Course in the History of Geometry in the 19th Century},
	author={Gray, Jeremy},
	year={2007},
	publisher={Springer},
	address={London},
	doi={10.1007/978-0-85729-060-1}
}

@book{cantor1932,
	title={Gesammelte Abhandlungen},
	author={Cantor, Georg},
	publisher={Verlag von Julius Springer},
	year={1932},
	address={Berlin}
}

@misc{wittmann2018,
	author       = {Wittmann, Erich Christian},
	title        = {{Das Wesen der Mathematik liegt in ihrer Freiheit: Fachliche Argumente f{\"u}r individuelle Lern- und L{\"o}sungswege}},
	year         = {2018},
	organization = {Technische Universit{\"a}t Dortmund, Projekt \glqq mathe 2000\grqq},
	url          = {https://wwwold.mathematik.tu-dortmund.de/ieem/mathe2000/pdf/Symp18/wittmann.pdf},
	note         = {Begleittext zum Symposium \glqq mathe 2000\grqq{} (Abgerufen am 2026-01-23)},
	langid       = {german}
}

@book{hoyrup2013,
	title={Lengths, Widths, Surfaces: a Portrait of Old Babylonian Algebra and its Kin},
	author={H{\o}yrup, Jens},
	year={2013},
	publisher={Springer Science \& Business Media},
	address={New York}
}

@book{friberg2016,
	title={New Mathematical Cuneiform texts},
	author={Friberg, J{\"o}ran and Al-Rawi, Farouk N.H.},
	year={2016},
	publisher={Springer},
	address={Cham},
	doi={10.1007/978-3-319-32219-5}
}

@book{kant_krv,
  author = {Immanuel Kant},
  title = {Kritik der reinen Vernunft},
  year = {1781/1787},
  note = {Zitiert nach der Originalausgabe A 290 / B 347 (Anhang zur Amphibolie der Reflexionsbegriffe)},
  publisher = {Suhrkamp},
  address = {Frankfurt am Main}
}

@book{Peano1889,
  title={Arithmetices principia: nova methodo},
  author={Peano, G.},
  year={1889},
  publisher={Fratres Bocca},
  address={Torino}
}

@book{peano1894,
  title={Formulaire de math{\'e}matiques},
  author={Peano, G.},
  year={1905},
  publisher={Fratres Bocca},
  address={Torino}
}

@book{ifrah2000,
	title={The universal history of numbers},
	author={Ifrah, Georges},
	year={2000},
	publisher={Harvill Press},
	address={London}
}

@incollection{Chapter17NoughtMatterstheHistoryandPhilosophyofZero,
	author = {Paul Ernest},
	title = {Chapter 17 Nought Matters: the History and Philosophy of Zero},
	booktitle = {The Origin and Significance of Zero: An Interdisciplinary Perspective},
	year = {2024},
	publisher = {Brill},
	address = {Leiden, Niederlande},
	pages= {306 --342}
}

@book{TheOriginandSignificanceofZero,
	author = {Peter Gobets and Robert Lawrence Kuhn},
	title = {The Origin and Significance of Zero: An Interdisciplinary Perspective},
	year = {2024},
	publisher = {Brill},
	address = {Leiden, Niederlande},
	isbn = {978-90-04-69156-8},
	doi = {10.1163/9789004691568}
}

@book{sharer2006,
	title={The ancient maya},
	author={Sharer, Robert J and Traxler, Loa P},
	year={2006},
	publisher={Stanford University Press},
	address={Stanford, CA}
}

@article{kanamori2003,
	title={The empty set, the singleton, and the ordered pair},
	author={Kanamori, Akihiro},
	journal={Bulletin of Symbolic Logic},
	volume={9},
	number={3},
	pages={273--298},
	year={2003},
	publisher={Cambridge University Press},
	doi={10.2178/bsl/1058448674}
}

@article{johnson1972,
	title={The genesis and development of set theory},
	author={Johnson, Phillip E},
	journal={The Two-Year College Mathematics Journal},
	volume={3},
	number={1},
	pages={55--62},
	year={1972},
	publisher={Taylor \& Francis},
	doi={10.2307/3026592}
}

@book{boole1854,
	title={An investigation of the laws of thought: on which are founded the mathematical theories of logic and probabilities},
	author={Boole, George},
	volume={2},
	year={1854},
	publisher={Walton and Maberly},
	address={London}
}

@article{fischbein1998,
	title={The mathematical concept of set and the'collection'model},
	author={Fischbein, Efraim and Baltsan, Madlen},
	journal={Educational studies in mathematics},
	volume={37},
	number={1},
	pages={1--22},
	year={1998},
	publisher={Springer},
	doi={10.1023/A:1003193213506}
}

@book{dedekind2012,
  title={Was sind und was sollen die Zahlen?},
  author={Dedekind, R.},
  isbn={9781108050388},
  series={Cambridge Library Collection - Mathematics},
  year={1887},
  publisher={Cambridge University Press},
  address={Cambridge},
  doi={10.1017/CBO9781139494390}
}

@book{Quine2013,
	title={Mengenlehre und ihre Logik},
	author={Quine, Willard Orman},
	volume={10},
	year={1973},
	publisher={Springer-Verlag},
	address={Berlin}
}

@article{nickel2007,
	title={{Mathematik: Schöpfung (fast) aus dem Nichts}},
	author={Nickel, Gregor},
	journal={Diagonal. Zeitschrift der Universit{\"a}t Siegen},
	pages={57--61},
	year={2007}
}

@book{pukelsheim2005,
  title={ Perspektivität und Unendlichkeit},
  author={Pukelsheim, Friedrich and Schwaetzer, Harald},
  publisher={Roderer Verlag},
  volume={29},
  year={2005}
}

@book{muller2010,
  title={Das Mathematikverst{\"a}ndnis des Nikolaus von Kues},
  author={Müller, Tom},
  publisher={Mathematische, Naturwissenschaftliche und Philosophisch-theologische Dimensionen. Mitteilungen und Forschungsbeitr{\"a}ge der Cusanus-Gesellschaft},
	address={Regensburg},
  year={2010}
}

@book{kant2013,
  title={Versuch den Begriff der negativen Grö\ss{}en in die Weltweisheit einzuführen},
  subtitle={Erstdruck: Königsberg (Kanter) 1763},
  author={Kant, Immanuel},
  publisher={Edition Holzinger},
    address={Berlin},
  year={2013}
}

@article{codd1979,
  author = {Codd, Edgar F.},
  title = {Extending the database relational model to capture more meaning},
  journal = {ACM Transactions on Database Systems (TODS)},
  volume = {4},
  number = {4},
  pages = {397--434},
  year = {1979},
  publisher = {ACM New York, NY, USA},
  doi = {10.1145/320107.320109}
}

@book{rathgeb2016,
  title={George Spencer Brown’s Laws of form between mathematics and philosophy : content - genesis - validitys},
  author={Rathgeb, Martin},
  year={2016},
  publisher={universi - Universitätsverlag Siegen},
  address={Siegen}
}

@incollection{nickel2018,
  author       = {Nickel, Gregor},
  title        = {{Kurzschlüsse oder fruchtbare wechselseitige Irritationen. Begegnungen von Mathematik und Theologie bei Nikolaus von Kues und Georg Cantor.}},
  editor       = {Hoff, Gregor Maria and Korber, Nikolaus},
  booktitle    = {Interdisziplin{\"a}re Forschung?: Ann{\"a}herungen an einen strapazierten Begriff},
  publisher    = {Verlag Karl Alber},
  year         = {2018},
  pages        = {150--187}
}

@incollection{nickel2021,
  author       = {Nickel, Gregor},
  title        = {Geist und {Zahl} (De mente c. 6)},
  editor       = {Mandrella, Isabelle},
  booktitle    = {Nicolaus Cusanus: Der Laie {\"u}ber den Geist/Idiota de mente},
  publisher    = {Walter de Gruyter GmbH \& Co KG},
  year         = {2021},
  pages        = {107--130}
}

@book{stwolfgang,
	title={Wolfgang und Johann Bolyai Geometrische Untersuchungen},
	author={St{\"a}ckel, P.},
	year={1913},
	publisher={Teubner},
	address={Leipzig}
}

@book{barna1970,
	title={A magyarorsz{\'a}gi matematika t{\"o}rt{\'e}nete [Geschichte der Mathematik in Ungarn]},
	author={Sz{\'e}n{\'a}ssy, Barna},
	publisher={Akad{\'e}miai Kiad{\'o}},
	address={Budapest},
	year={1970}
}

@book{2007platon,
  author    = {Platon},
  title     = {Sophistes. Hrsg. von U. Wolf, F. Schleiermacher und C. Iber},
  isbn      = {9783518270042},
  series    = {Suhrkamp Studienbibliothek},
  publisher = {Suhrkamp},
  year      = {2007},
  address   = {Frankfurt am Main}
}

@book{rotman1987,
  title={Signifying Nothing: The Semiotics of Zero},
  author={Rotman, Brian},
  year={1987},
  publisher={Stanford University Press},
  address={Stanford, CA}
}

\end{document}